\documentclass[12pt]{amsproc}
\theoremstyle{plain}
\newtheorem{theorem}{Theorem}[section]
\newtheorem{corollary}[theorem]{Corollary}
\newtheorem{example}[theorem]{Example}

\newtheorem{proposition}[theorem]{Proposition}
\newtheorem{lemma}[theorem]{Lemma}

\theoremstyle{definition}
\newtheorem{definition}[theorem]{Definition}

\theoremstyle{remark}

\begin{document}

\title{A module frame concept for Hilbert C*-modules}
\author[M.~Frank]{Michael Frank}
\address{Universit\"at Leipzig, Mathematisches Institut, D-04109 Leipzig, F.R. Germany}
\email{frank@mathematik.uni-leipzig.de}
\thanks{Both authors were supported by grants from the NSF}
\thanks{To appear in: {\it Functional and Harmonic Analysis of Wavelets (San
Antonio, TX, Jan.~1999)}, Contemp.~Math. {\bf 247}, 207-233,
A.M.S., Providence, RI, 2000.}
\author[D.~R.~Larson]{David R.~Larson}
\address{Dept.~Mathematics, Texas A{\&}M Univ., College Station, TX 77843, U.S.A.}
\email{larson@math.tamu.edu}
\begin{abstract}
The goal of the present paper is a short introduction to a general mo\-dule
frame theory in C*-algebras and Hilbert C*-modules. The reported investigations
rely on the idea of geometric dilation to standard Hilbert C*-modules over
unital C*-algebras that possess orthonormal Hilbert bases, and of reconstruction
of the frames by projections and other bounded module operators with suitable
ranges. We obtain frame representation and decomposition theorems, as well as
similarity and equivalence results. The relative position of two and more frames
in terms of being complementary or disjoint is investigated in some detail.
In the last section some recent results of P.~G.~Casazza are generalized to
our setting. The Hilbert space situation appears as a special case. For the
details of most of the proofs we refer to our basic publication [8].
\end{abstract}
\maketitle

Frames serve as a replacement for bases in Hilbert spaces that guarantee canonical
reconstruction of every element of the Hilbert space by the reconstruction
formula, however, giving up linear independence of the elements of the
generating frame sequence. They appear naturally as wavelet generated and
Weyl-Heisenberg / Gabor frames since often sequences of this type do not
become orthonormal or Riesz bases, \cite{HanLarson,Cas/Chr,Hei/Wal}.
Similarly, the concept of module frames has become a fruitful tool in
C*-algebra theory, e.~g.~for the description of conditional expectations of
finite (Jones) index \cite{Wata,PP,BDH,KiFr98} or for the investigation of
Cuntz-Krieger-Pimsner algebras \cite{DPZ,KPW}. In this part of the literature
they are called 'quasi-bases' or 'bases'.
Inspired by the results by Deguang Han and D.~R.~Larson in \cite{HanLarson}
and by investigations by M.~A.~Rieffel \cite{Rie2} the two authors
investigated completions of algebraic tensor products $A \odot H$ of a
C*-algebra $A$ and a (separable) Hilbert space $H$ that arise from the
extension of the scalar product on $H$ to an $A$-valued inner product on
$A \odot H$ as well as their direct summands for the existence of module
frames. Surprisingly, a general concept of module frames for this class of
Hilbert $A$-modules was obtained and most of the theorems valid for the
Hilbert space situation are extendable to theorems for Hilbert C*-modules
undergoing only minor additional changes, \cite{LaFr98}. The results
are useful in wavelet theory, in C*-theory and (in their commutative
reinterpretation) in vector and Hilbert bundle theory \cite{Fr99}.

\smallskip
The goal of the present paper is to give an introduction to the more general
concept of module frames and to indicate further results extending those
explained in \cite{LaFr98}. We start with a brief description of M.~A.~Rieffel's
link of C*-theory to wavelet theory via group representations
giving the motivation for our investigations. In section two we give an
introduction to the basic notions of Hilbert C*-module theory used in the
sequel. Sections three and four are devoted to the crucial results showing
the wide field of analogies between frame theory in Hilbert spaces and in
Hilbert C*-modules and the few obstacles that arise in the generalized
situations. A detailed account to inner sums of frames and their properties,
including considerations about existence and uniqueness, is the goal of
section five and six. There we reproduce Hilbert space results of
\cite{HanLarson} in our more general setting with slight modifications.
We close our explanations with investigations on module frame
representations and decompositions generalizing results by P.~G.~Casazza
\cite{Casazza:1} that describe the Hilbert space situation.

\section{Module frames and multi-resolution analysis wavelets}

Following M.~A.~Rieffel \cite{Rie2}
assume the situation of a wavelet sequence generated by a multi-resolution
analysis in a Hilbert space $L_2({\mathbb R}^n)$. Denote the mother wavelet
by $\phi \in L_2({\mathbb R}^n)$, $\| \phi \|_2 =1$, and consider
${\mathbb R}^n$ as an additive group. The second group appearing in the picture
is $\Gamma = {\mathbb Z}^n$ acting on $L_2({\mathbb R}^n)$ by translations in
the domains of functions, i.e.~mapping $\phi(x)$ to $\phi(x-p)$ for $x \in
{\mathbb R}^n$ and $p \in {\mathbb Z}^n$. By the origin of the mother wavelet
$\phi$ all the different ${\mathbb Z}^n$-translates of $\phi$ are pairwise
orthogonal,
  \[
    \int_{{\mathbb R}^n} \overline{\phi(x-q)}\phi(x-p) \, dx = \delta_{qp} \, .
  \]
Introducing the group C*-algebras $A=C^*({\mathbb Z}^n)$ of the additive
discrete group ${\mathbb Z}^n$ into the picture and interpreting the set of
all ${\mathbb Z}^n$-translates of $\phi$ as elements of the $*$-algebra
$C_c({\mathbb R}^n)$ we obtain a right action of $A$ on $C_c({\mathbb R}^n)$
by convolution and an $A$-valued inner product there defined by
  \[
    \langle \phi,\psi \rangle_A(p) := \int_{{\mathbb R}^n} \overline{\phi(x)}
    \psi(x-p) \, dx
  \]
for $\phi, \psi \in C_c({\mathbb R}^n)$ and $p \in {\mathbb Z}^n$,
see section two for a detailed explanation of the notion. The
completion of $C_c({\mathbb R}^n)$ with respect to the norm $\| \phi \| :=
\| \langle \phi,\phi \rangle_A \|_A^{1/2}$ is a right Hilbert C*-module
$\mathcal H = \overline{C_c({\mathbb R}^n)}$ over $A$.

Considering the dual Fourier transformed picture things become mathematically
easier. The C*-algebra $A=C^*({\mathbb Z}^n)$ is transformed to the C*-algebra
$B=C({\mathbb T}^n)$ of continuous functions on the $n$-torus. The right
action of $A$ on $\mathcal H$ by convolution becomes a right action of $B$ on
$\mathcal H$ by pointwise multiplication. Moreover, $\mathcal H = \overline{C_c
({\mathbb R}^n)}$ coincides with the set $B \phi$, i.e.~it is a singly generated
free $B$-module with $B$-valued inner product
  \[
    \langle \phi, \psi \rangle_B(t) := \sum_{p \in {\mathbb Z}^n}
                   (\overline{\phi}\psi)(t-p)
  \]
for $t \in {\mathbb R}^n$. The set $\{ \phi \}$ consisting of one element is
a module frame, even a module Riesz basis. However, for $n \geq 2$ there
exist non-free $B$-modules that are direct orthogonal summands of $\mathcal H
= B$, cf.~\cite{Rie2} for their construction. For them module bases might not
exist and module frames consist of more than one element.

In a similar manner multi-wavelets give rise to Hilbert $B$-modules $B^k$ of
all $k$-tuples with entries from $B$ and coordinate-wise operations. Since
they are also free finitely generated $B$-modules the pairwise orthogonal
generating elements of the multi-wavelet $\{ \phi_1, ...,\phi_k \}$ form a
module frame for $B^k$. Considering infinite sets of pairwise orthogonal
generating elements in the context of a multi-resolution analysis we have to
find suitable closures of the algebraic tensor product $A \odot l_2$ to form
well-behaved Hilbert C*-modules. Since norm-convergence and weak convergence
are in general different concepts in an infinite-dimensional C*-algebra $B$
(whereas both they coincide in $\mathbb C$) some more investigations have to
be carried out.

\section{A brief introduction to Hilbert C*-modules}

Hilbert C*-modules arose as generalizations of the notion 'Hilbert space'.
The basic idea was to consider modules over C*-algebras instead of linear spaces
and to allow the inner product to take values in the C*-algebra of coefficients
being C*-(anti-)linear in its arguments. The structure was first used by
I.~Kaplansky \cite{Kapl} in 1952 and more carefully investigated by
M.~A.~Rieffel \cite{Rie} and W.~L.~Paschke \cite{Pa1} later in 1972/73.
For comprehensive accounts we refer to the lecture notes of E.~C.~Lance
\cite{La} and to the book of N.-E.~Wegge-Olsen \cite[ch.~15, 16]{NEWO}.
We give only a brief introduction to the theory of Hilbert C*-modules to make
our explanations self-contained. Throughout the present paper we assume the
C*-algebra of coefficients of Hilbert C*-modules to be unital.

\begin{definition}
 Let $A$ be a (unital) C*-algebra and $\mathcal H$ be a (left) $A$-module.
 Suppose that the linear structures given on $A$ and $\mathcal H$ are
 compatible, i.e. $\lambda (a x) =(\lambda a)x = a(\lambda x)$ for every
 $\lambda \in {\mathbb C}$, $a \in A$ and $x \in \mathcal H$. If there
 exists a mapping $\langle .,. \rangle: {\mathcal H} \times {\mathcal H}
 \rightarrow A$ with the properties

   \newcounter{marke}
   \begin{list}{(\roman{marke})}{\usecounter{marke}}
      \item $\langle x,x \rangle \geq 0$ for every $x \in \mathcal H$,
      \item $\langle x,x \rangle =0$ if and only if $x=0$,
      \item $\langle x,y \rangle = \langle y,x \rangle^*$ for every $x,y \in
            \mathcal H$,
      \item $\langle ax,y \rangle = a \langle x,y \rangle$ for every $a \in A$,
            every $ x,y \in \mathcal H$,
      \item $\langle x+y,z \rangle = \langle x,z \rangle + \langle y,z \rangle$
            for every $x,y,z \in \mathcal H$.
   \end{list}

 Then the pair $\{ \mathcal H, \langle .,. \rangle \}$ is called a
 {\it (left) pre-Hilbert $A$-module}. The map $\langle .,. \rangle$ is said to
 be an {\it $A$-valued inner product}.
 If the pre-Hilbert $A$-module $\{ \mathcal H, \langle .,. \rangle \}$ is complete
 with respect to the norm $\|x\| = \| \langle x,x \rangle \|^{1/2}$ then
 it is called a {\it Hilbert $A$-module}.

 In case $A$ is unital the Hilbert $A$-module $\mathcal H$ is {\it
 (algebraically) finitely generated} if there exists a finite set $\{ x_i
 \}_{i \in \mathbb N} \subset \mathcal H$ such that $x = \sum_i a_i x_i$ for
 every $x \in \mathcal H$ and some coefficients $\{ a_i \} \subset A$.
 If $A$ is unital the Hilbert $A$-module $\mathcal H$ is {\it countably
 generated} if there exists a countable set $\{ x_i \}_{i \in \mathbb N}
 \subset \mathcal H$ such that the set of all finite $A$-linear combinations
 $\{ \sum _j a_jx_j \}$, $\{ a_i \} \subset A$, is norm-dense in $\mathcal H$.
\end{definition}

\begin{example} {\rm
  Denote by ${\mathcal H} =A^n$ the set of all $n$-tuples with entries of $A$,
  where the addition is the position-wise addition derived from $A$, the action
  of $A$ on $\mathcal H$ is the multiplication of every entry by a fixed
  element of $A$ from the left and the $A$-valued inner product is defined by
  the formula
   \[
       \langle a,b \rangle = \sum_{i=1}^n a_i b_i^*
   \]
  for $a=(a_1, ...,a_n), b=(b_1, ... ,b_n) \in \mathcal H$. This Hilbert
  $A$-module is a free finitely generated one with a standard orthonormal
  basis.   }
\end{example}

This kind of examples plays a crucial rule in Hilbert C*-module theory. It
allows to characterize finitely generated Hilbert C*-modules precisely as the
projective finitely generated C*-modules.

\begin{theorem} \label{th1}
  Let $A$ be a unital C*-algebra.
  Every algebraically finitely generated Hilbert $A$-module $\mathcal H$ is
  an orthogonal summand of some free Hilbert $A$-module $A^n$ for a finite
  number $n$.
\end{theorem}

\begin{example}  {\rm
  Let $H$ be a Hilbert space, $A \odot H$ be the algebraic tensor product
  of $H$ and the C*-algebra $A$ of coefficients, and define the $A$-valued
  inner product by
   \[
      \langle a \otimes h, b \otimes g \rangle = a \langle h,g \rangle_H b^*
   \]
  for $a,b  \in A$, $h,g \in H$. Then $A \odot H$ becomes a pre-Hilbert
  $A$-module.

  \noindent
  Naturally, $A^n \cong A \odot {\mathbb C}^n$ for any $n \in \mathbb N$.
  Moreover, set $l_2(A)$ to be the norm-completed algebraic tensor product
  $\overline{A \odot l_2}$ that possesses an alternative description as
   \[
      l_2(A) = \{ a=\{a_i\}_{i \in {\mathbb N}} \, : \,\, {\rm the} \,\,
      {\rm sum} \,\,\, {\sum}_{j=1}^{\infty} a_j a_j^*
      \,\,\, {\rm converges} \,\, {\rm w.r.t.} \,\, \|.\|_A \}
   \]
  with $A$-valued inner product $\langle a,a \rangle = \sum_j a_j a_j^*$. }
\end{example}

The Hilbert $A$-module $l_2(A)$ serves as an universal environment for
countably gene\-rated Hilbert $A$-modules that can be described as (special)
orthogonal summands. This fact was first observed by G.~G.~Kasparov
\cite{Kas} in 1980.

\begin{theorem} \label{th2} {\rm (Stabilization theorem)}
  Let $A$ be a unital C*-alge\-bra.
  Every countably generated Hilbert $A$-module $\mathcal H$ possesses an
  embedding into $l_2(A)$ as an orthogonal summand in such a way that the
  orthogonal complement of it is isometrically isomorphic to $l_2(A)$ again,
  i.e.~$\mathcal H \oplus l_2(A) = l_2(A)$.
\end{theorem}

The Hilbert $A$-modules $A^n$, $n \in \mathbb N$, and $l_2(A)$ possess
canonical orthonormal bases. However, not every Hilbert C*-module has an
orthonormal basis (cf.~Example \ref{example1}), and some replacement for this
notion would be helpful to describe the inner structure of Hilbert C*-modules.

\section{Frames in Hilbert C*-modules}

The concept that replaces the notion of a bases and strengthens the simple
notion of a generating set is that one of a module frame. Of course, the
existence of module frames has to be shown separately.

\begin{definition}
  Let $A$ be a unital C*-algebra and ${\mathbb J}$ be a finite or countable
  index set.
  A sequence $\{ x_j : j \in {\mathbb J} \}$ of elements in a Hilbert
  $A$-module $\mathcal H$ is said to be a {\it frame} if there are real
  constants $C,D > 0$ such that
  \begin{equation} \label{ineq-frame}
      C \cdot \langle x,x \rangle \leq
      \sum_{j} \langle x,x_j \rangle \langle x_j,x \rangle \leq
      D \cdot \langle x,x \rangle
  \end{equation}
  for every $x \in \mathcal H$. The optimal constants (i.e.~maximal for $C$
  and minimal for $D$) are called {\it frame bounds}. The frame $\{ x_j :
  j \in {\mathbb J} \}$ is said to be a {\it tight frame} if $C=D$, and said
  to be {\it normalized} if $C=D=1$.
  We consider {\it standard} (normalized tight) frames in the main for which
  the sum in the middle of the inequality (\ref{ineq-frame}) always converges
  in norm. For non-standard frames the sum in the middle converges only weakly
  for at least one element of $\mathcal H$.

  \smallskip
  A sequence $\{ x_j \}_j$ is said to be a {\it standard Riesz basis of
  $\mathcal H$} if it is a standard frame and a generating set with the
  additional property that $A$-linear combinations $\sum_{j \in S} a_jx_j $
  with coefficients $\{ a_j \} \in A$ and $S \in {\mathbb J}$ are equal to zero
  if and only if in particular every summand $a_jx_j$ equals zero for $j \in S$.
  A generating sequence $\{ x_j \}_j$ with the described additional property
  alone is called a {\it Hilbert basis of $\mathcal H$}.

  An {\it inner summand of a standard Riesz basis} of a Hilbert $A$-module
  $\mathcal L$ is a sequence $\{ x_j \}_j$ in a Hilbert $A$-module $\mathcal H$
  for which there is a second sequence $\{ y_j \}_j$ in a Hilbert $A$-module
  $\mathcal K$ such that ${\mathcal L} \cong {\mathcal H} \oplus {\mathcal K}$
  and the sequence consisting of the pairwise orthogonal sums $\{ x_j
  \oplus y_j \}_j$ in the Hilbert $A$-module ${\mathcal H} \oplus {\mathcal K}$
  is the initial standard Riesz basis of $\mathcal L$.

  \smallskip
  Two frames $\{ x_j \}_j$, $\{ y_j \}_j$ of Hilbert $A$-modules $H_1$, $H_2$,
  respectively, are {\it unitarily equivalent} (resp., {\it similar}) if there
  exists a unitary (resp., invertible adjointable) linear operator $T: H_1 \to
  H_2$ such that $T(x_j)=y_j$ for every $j \in \mathbb J$.
\end{definition}

Recalling the standard identifications $A^n \cong A \odot {\mathbb C}^n$ and
$l_2(A) \cong \overline{A \odot l_2}$ we observe that for every (normalized
tight) frame $\{ x_j \}_j$ of a Hilbert space $H$ the sequence $\{ 1_A \otimes
x_j \}_j$ is a standard (normalized tight) module frame of the
Hilbert $A$-module $\mathcal H = \overline{A \odot H}$ with the same frame
bounds. So standard module frames exist in abundance in the canonical
Hilbert $A$-modules. To show the existence of standard module frames in
arbitrary finitely or countably generated Hilbert $A$-modules we have the
following fact:

\begin{theorem}   \label{th3}
  For every $A$-linear partial isometry $V$ on $A^n$ (or $l_2(A)$) the image
  sequence $\{ V(e_j) \}_j$ of the standard orthonormal basis $\{ e_j \}_j$
  is a standard normalized tight frame of the image $V(A^n)$ (or $V(l_2(A))$).
  Consequently, every algebraically finitely generated (or countably generated)
  Hilbert $A$-module $\mathcal H$ possesses a standard normalized tight frame.
\end{theorem}

The verification of this statement is straightforward starting from a standard
decomposition of elements with respect to some orthonormal basis and applying
the partial isometry to both the sides of the equality. A reference to the
Theorems \ref{th1} and \ref{th2} gives the existence of frames since
projections are special partial isometries.

We point out that a standard Riesz basis $\{ x_j \}_j$ which is a normalized
tight frame is automatically an orthogonal Hilbert basis with the property
$\langle x_j,x_j \rangle = \langle x_j,x_j \rangle^2$ for any $j \in \mathbb
J$, cf.~\cite[Prop.~2.2]{LaFr98}.
Finally, we show that certain countably generated Hilbert C*-modules do not
possess any standard (orthogonal) Riesz basis.

\begin{example}  \label{example1} {\rm
Take the C*-algebra $A=C([0,1])$
of all continuous functions on the unit interval $[0,1]$ that is equipped with
the usual topology. Consider the Hilbert $A$-module $\mathcal H = C_0((0,1])$
together with the $A$-valued inner product  $\langle f,g \rangle (t) =
f(t)\overline{g(t)}$, $t \in [0,1]$. The set $\{ f(t)=t \} \in \mathcal H$ is
a topologically (but not algebraically) generating set and naturally
orthonormal. However, it is not a frame, since the lower frame bound $C$
equals zero. The latter is true for any set consisting of exactly one generator.
If we take generating sets consisting of two (or more) generators then the
requirement to them to be orthogonal translates into the fact that there
must exist points $t_0 \in (0,1]$ where all generators vanish. This is a
contradiction to the first property of the set to be generating. Conversely,
if the set is generating and contains more than one element, then at least
two of the generators possess a support subinterval in common. This allows a
non-trivial representation of the zero function, so these sets cannot be
bases. At this point we have exhausted all possibilities to build a Riesz
basis for the Hilbert $A$-module $\mathcal H$. However, the Hilbert C*-module
under considerations possesses a standard normalized tight frame, see Example
\ref{example2}.   }
\end{example}

\section{Reconstruction formulae, frame transform and frame operator}

One of the main properties of normalized tight frames (standard or
non-standard) is the validity of a canonical reconstruction formula for all
elements of the spanned Hilbert C*-module. Conversely, every generating
sequence of a countably generated Hilbert C*-module satisfying this type
reconstruction formula is naturally a normalized tight frame for it. For
more general frames the situation is slightly more complicated since a dual
frame is needed to build up a reconstruction formula. Its existence can only
be guaranteed if the initial frame is supposed to be standard. Since the
methods of the proofs of the theorems below involve a lot of technicalities
and deeper knowledge of Hilbert C*-module theory we refer to \cite{LaFr98}
for a detailed theoretical approach. The Hilbert space situation is investigated
in great detail by Deguang Han and D.~R.~Larson in \cite{HanLarson}.

\begin{theorem}  {\rm (\cite[Th.~4.1]{LaFr98} )} \label{th4}
  Let $A$ be a unital C*-algebra, $\mathcal H$ be a finitely or countably
  generated Hilbert $A$-module and $\{ x_j \}_j$ be a normalized tight frame
  of $\mathcal H$. Then the {\it reconstruction formula}
   \begin{equation}  \label{eq1}
     x= \sum_j \langle x,x_j \rangle x_j
   \end{equation}
  holds for every $x \in \mathcal H$ in the sense of convergence w.r.t.~the
  topology that is induced by the set of semi-norms
  $\{ |f(\langle .,. \rangle )|^{1/2} \, : \, f \in A^* \}$.
  The sum converges always in norm if and only if the frame $\{ x_j \}_j$ is
  standard.      \newline
  Conversely, a finite set or a sequence $\{ x_j \}_j$ satisfying the formula
  (\ref{eq1}) for every $x \in \mathcal H$ is a normalized tight frame of
  $\mathcal H$.
\end{theorem}

\begin{example} \label{example2} {\rm
  Consider the Hilbert $C([0,1])$-module $C_0((0,1])$ and the standard
  $A$-valued inner product derived from multiplication and involution in $A$
  again, cf.~Example \ref{example1}. Take the functions
  \begin{equation*}
     f_1(t)  = \left\{ \begin{array}{ccl}
               0 & : &  t \in [0,1/2] \\
               g(t-1/2) & : & t \in [1/2,1]
               \end{array} \right.    \, ,
  \end{equation*}
  \begin{equation*}
      f_n(t)  = \left\{ \begin{array}{ccl}
               0 & : & t \in [0,1/2^n] \\
               g(2^{(n-1)}t -1/2) & : & t \in [1/2^n,1/2^{(n-1)}] \\
               g(2^{(n-2)}t) & : & t \in [1/2^{(n-1)},1/2^{(n-2)}]
               \end{array} \right. \, ,
  \end{equation*}
  where $g(t)=2t$ for $t \in [0,1/2]$ and $g(t)=2 \sqrt{t-t^2}$ for $t \in
  [1/2,1]$. The sequence $\{ f_i \}_i$ forms a standard normalized tight frame
  for the Hilbert $C([0,1])$-module $C_0((0,1])$. This can be shown best
  verifying the reconstruction formula and checking their convergence in norm
  for any concrete element of the Hilbert C*-module.   }
\end{example}

With some experience in Hilbert C*-module theory the following crucial
fact is surprising because of the generality in which it holds. The existence
and the very good properties of the frame transform of standard frames give
the chance to get far reaching results analogous to those in the Hilbert space
situation. Again, the proof is more complicated than the known one in the
classical Hilbert space case, cf.~\cite{HanLarson}.

\begin{theorem}   {\rm (\cite[Th.~4.2]{LaFr98} )}   \label{th5}
  Let $A$ be a unital C*-algebra, $\mathcal H$ be a finitely or countably
  generated Hilbert $A$-module and $\{ x_j \}_j$ be a standard frame
  of $\mathcal H$. The {\it frame transform of the frame $\{ x_j \}_j$} is
  defined to be the map
    \[
       \theta: \mathcal H \rightarrow l_2(A) \quad , \qquad \theta(x)
       = \{ \langle x,x_j \rangle \}_j
    \]
  that is bounded, $A$-linear, adjointable and fulfills $\theta^*(e_j)=x_j$
  for a standard orthonormal basis $\{ e_j \}_j$ of the Hilbert $A$-module
  $l_2(A)$ and all $j \in \mathbb J$.   \newline
  Moreover, the image $\theta(\mathcal H)$ is an orthogonal summand of
  $l_2(A)$. For normalized tight frames we additionally get $P(e_j)=\theta
  (x_j)$ for any $j \in \mathbb J$, and $\theta$ is an isometry in that case.
\end{theorem}

The frame transform $\theta$ is the proper tool for the description of
standard frames.

\begin{theorem}  {\rm (\cite[Th.~6.1]{LaFr98} )} \label{th6}
  Let $A$ be a unital C*-algebra, $\mathcal H$ be a finitely or countably
  generated Hilbert $A$-module and $\{ x_j \}_j$ be a standard frame
  of $\mathcal H$. Then the {\it reconstruction formula}
    \[
         x= \sum_j \langle x,S(x_j) \rangle x_j
    \]
  holds for every $x \in \mathcal H$ in the sense of norm convergence, where
  $S:=(\theta^*\theta)^{-1}$ is positive and invertible.
\end{theorem}

The sequence $\{ S(x_j) \}_j$ is a standard frame again, the {\it canonical
dual frame of} $\{ x_j \}_j$. The operator $S$ is called the {\it frame
operator of $\{ x_j \}_j$} on $\mathcal H$. The Theorems  \ref{th4} and
\ref{th6} bring to light some key properties of frame sequences:

\begin{corollary} {\rm (\cite[Th.~5.4, 5.5]{LaFr98}) }
  Every standard frame of a finitely or countably generated Hilbert $A$-module
  is a set of generators.  \newline
  Every finite set of algebraic generators of a finitely generated Hilbert
  $A$-module is a (standard) frame.
\end{corollary}

Sometimes the reconstruction formula of standard frames is valid with other
(standard) frames $\{ y_i \}_i$ instead of $\{ S(x_i) \}_i$. They are said to
be {\it alternative dual frames of $\{ x_i \}_i$}. However, the canonical
dual frame is in some sense optimal:

\begin{corollary} {\rm (\cite[Prop.~6.7]{LaFr98}) }
    If $x=\sum_j \langle x,y_j \rangle x_j$ for every $x \in \mathcal H$ and
    a different frame $\{ y_j \}_j \subset \mathcal H$, then we have the
    inequality
      \[
           \sum_j \langle x,S(x_j) \rangle \langle S(x_j),x \rangle <
           \sum_j \langle x,y_j \rangle \langle y_j,x \rangle < \infty \, .
      \]
    that is valid for every $x \in \mathcal H$ in the positive cone of $A$.
\end{corollary}

We close our considerations on the frame transform arising from standard
frames with a statement on the relation between unitary equivalence (or
similarity) of frames and the characteristics of the image of the frame
transform $\theta$.

\begin{theorem} {\rm (\cite[Th.~7.2]{LaFr98}) } \label{th4.2}
  Let $A$ be a unital C*-algebra and $\{ x_j \}_j$ and $\{ y_j \}_j$ be
  standard (normalized tight) frames of Hilbert $A$-modules ${\mathcal H}_1$
  and ${\mathcal H}_2$, respectively. The following conditions are
  equivalent:

  \newcounter{cou008}
  \begin{list}{(\roman{cou008})}{\usecounter{cou008}}
   \item The frames $\{ x_j \}_j$ and $\{ y_j \}_j$ are unitarily equivalent
         or similar.
   \item Their frame transforms $\theta_1$ and $\theta_2$ have the same range
         in $l_2(A)$.
   \item The sums $\sum_j a_j x_j$ and $\sum_j a_j y_j$ equal zero for exactly
         the same Banach A-submodule of sequences $\{ a_j \}_j$ of $l_2(A)$.
  \end{list}
\end{theorem}

\section{Complementary frames and inner sums of frames}

The aim of the present section is a description of existing module frame
complements that allow to understand standard frames as parts of orthonormal
and Riesz bases of extended Hilbert C*-modules, sometimes up to unitary
equivalence. Throughout we use a special composition of two frames of
different Hilbert C*-modules that consists of pairs of frame elements strictly
with identical indices, the inner sum of them. Thinking about two orthonormal
bases of two different Hilbert spaces, as a particular case, this procedure
will not yield an orthonormal basis of the direct sum of these Hilbert spaces
at all. For frames, however, this may happen to be a frame for the direct sum,
again. The usefulness of this concept was first observed in \cite{HanLarson},
where the Hilbert space situation was worked out.

We start with standard normalized tight frames of Hilbert C*-modules. For
them some uniqueness results still hold, an important fact for further
investigations. Of course, because of the great variety of non-isomorphic
Hilbert C*-modules with the same cardinality of their sets of generators
we can state uniqueness only for appropriately fixed resulting direct sum
Hilbert C*-modules. (See \cite[5.1-5.7]{LaFr98} for details.)

\begin{proposition} \label{prop-existence} {\rm (Existence)}
  Let $A$ be a unital C*-algebra, $\mathcal H$ be a finitely (resp., countably)
  generated Hilbert $A$-module and $\{ x_j : j \in {\mathbb J} \}$ be a standard
  normalized tight frame of $\mathcal H$.
  Then there exists another countably generated Hilbert $A$-module $\mathcal M$
  and a standard normalized tight frame $\{ y_j : j \in {\mathbb J} \}$ in it
  such that
  \[
     \{ x_j \oplus y_j : j \in {\mathbb J} \}
  \]
  is an orthogonal Hilbert basis for the countably generated Hilbert
  $A$-module ${\mathcal H} \oplus {\mathcal M}$ with $\langle x_j \oplus y_j,
  x_j \oplus y_j \rangle = \langle x_j \oplus y_j, x_j \oplus y_j \rangle^2$.
  $\, \mathcal M$ can be chosen in such a way that ${\mathcal H} \oplus
  {\mathcal M} = l_2(A)$ and hence, $1_A = $  \linebreak[4]
  $=\langle x_j \oplus y_j, x_j
  \oplus y_j \rangle = \langle x_j \oplus y_j, x_j \oplus y_j \rangle^2$ for any
  $j \in \mathbb J$.

  If $\mathcal H$ is finitely generated and the index set ${\mathbb J}$ is finite
  then $\mathcal M$ can be chosen to be finitely generated, too, and ${\mathcal H}
  \oplus {\mathcal M} = A^N$ for $N = |{\mathbb J}|$.

  If $\{ x_j \}_j$ is already an orthonormal basis then ${\mathcal M}= \{ 0 \}$,
  i.e.~no addition to the frame is needed. If ${\mathbb J}$ is finite and
  $\mathcal M$ is not finitely generated (i.e.~we try to dilate to an infinitely
  generated Hilbert $A$-module) then infinitely many times $0_{\mathcal H}$ has
  to be added to the frame $\{ x_j \}_j$.
\end{proposition}

\begin{proposition} \label{prop-unique} {\rm (Uniqueness)}
  Let $A$ be a unital C*-algebra, $\mathcal H$ be a countably generated Hilbert
  $A$-module and $\{ x_j : j \in {\mathbb J} \}$ be a standard normalized tight
  frame for $\mathcal H$, where the index set ${\mathbb J}$ is countable or finite.
  Suppose, there exist two countably generated Hilbert $A$-modules $\mathcal M$,
  $\mathcal N$ and two normalized tight frames $\{ y_j : j \in {\mathbb J} \}$,
  $\{ z_j : j \in {\mathbb J} \}$ for them, respectively, such that
  \[
     \{ x_j \oplus y_j : j \in {\mathbb J} \} \, ,
     \, \{ x_j \oplus z_j : j \in {\mathbb J} \}
  \]
  are orthogonal Hilbert bases for the countably generated Hilbert $A$-modules
  ${\mathcal H} \oplus {\mathcal M}$, ${\mathcal H} \oplus {\mathcal N}$,
  respectively, where we have the value properties
   \begin{eqnarray*}
    \langle x_j \oplus y_j, x_j \oplus y_j \rangle = \langle x_j \oplus y_j,
    x_j \oplus y_j \rangle^2 \, & {\it and} & \\
    \langle x_j \oplus z_j, x_j \oplus z_j \rangle = \langle x_j \oplus z_j,
    x_j \oplus z_j \rangle^2 & . &
   \end{eqnarray*}
  If $\langle y_j,y_j \rangle_{{\mathcal M}} = \langle z_j,z_j
  \rangle_{{\mathcal N}}$ for every $j \in {\mathbb J}$, then there exists a
  unitary transformation
  $U: {\mathcal H} \oplus {\mathcal M} \to {\mathcal H}
  \oplus {\mathcal N}$ acting identically on $\mathcal H$, mapping $\mathcal
  M$ onto $\mathcal N$ and satisfying $U(y_j)=z_j$ for every $j \in {\mathbb J}$.

  The additional remarks of Proposition \ref{prop-existence} apply in the
  situation of finitely generated Hilbert $A$-modules appropriately.
\end{proposition}

Considering arbitrary standard frames we obtain similar results, however,
uniqueness cannot be established without further restrictions.

\begin{proposition}
  Let $\{ x_j : j \in {\mathbb J} \}$ be a standard frame of a finitely or
  countably generated Hilbert $A$-module $\mathcal H$. There exists a Hilbert
  $A$-module $\mathcal M$ and a normalized tight frame $\{y_j : j \in {\mathbb
  J} \}$ in $\mathcal M$ such that the sequence $\{ x_j \oplus y_j : j \in
  {\mathbb J} \}$ is a standard Riesz basis in ${\mathcal H} \oplus {\mathcal
  M}$ with the same frame bounds for both $\{ x_j \}$ and $\{ x_j \oplus y_j \}$.
  The Hilbert $A$-module $\mathcal M$ can be chosen in such a way that
  ${\mathcal H} \oplus {\mathcal M} = l_2(A)$. If $\mathcal H$ is finitely
  generated and the index set ${\mathbb J}$ is finite, then $\mathcal M$ can be
  chosen to be finitely generated, too, and ${\mathcal H} \oplus {\mathcal M}
  = A^N$ for $N=|{\mathbb J}|$.

  In general, $\mathcal M$ cannot be selected as a submodule of $\mathcal H$,
  and the resulting standard Riesz basis may be non-orthogonal. A uniqueness
  result like that one for standard normalized tight frames fails in general to
  be true.
\end{proposition}

Adding the basis property to the frame property, i.e.~taking standard Riesz
bases, we reobtain a classical characterization of them that is sometimes
taken as an alternate definition of standard Riesz bases.

\begin{proposition}
  Let $\{ x_j : j \in {\mathbb J} \}$ be a standard Riesz basis of a finitely or
  countably generated Hilbert $A$-module $\mathcal H$. Then it is the image of a
  standard normalized tight frame and Hilbert basis   \linebreak[4]
  $\{ y_j : j \in {\mathbb J}
  \}$ of $\mathcal H$ under an invertible adjointable bounded $A$-linear operator
  $T$ on ${\mathcal H}$, i.e.~of an orthogonal Hilbert basis $\{ y_j : j \in {\mathbb
  J} \}$ for which $\langle y_j,y_j \rangle = \langle y_j,y_j \rangle^2$
  holds for every $j \in \mathbb J$.

  Conversely, the image of a standard normalized tight frame and Hilbert basis
  $\{ y_j : j \in {\mathbb J} \}$ of $\mathcal H$ under an invertible adjointable
  bounded $A$-linear operator $T$ on ${\mathcal H}$ is a standard Riesz basis of
  $\mathcal H$.

  If a Hilbert $A$-module $\mathcal H$ contains a standard Riesz basis then it
  contains an orthogonal Hilbert basis $\{ x_j\}_j$ with the property that
  \linebreak[4]
  $x = \sum_j \langle x,x_j \rangle x_j$ for every element $x \in \mathcal H$.
\end{proposition}

We can summarize the results in this section in a short statement.

\begin{theorem}
  Standard frames are precisely the inner direct summands of standard Riesz
  bases of $A^n$, $n \in \mathbb N$, or $l_2(A)$. Standard normalized tight
  frames are precisely the inner direct summands of orthonormal Hilbert bases
  of $A^n$, $n \in \mathbb N$, or $l_2(A)$.
\end{theorem}

\section{Different types of disjointness of composite frames in inner sums}

The aim of this section is to consider sequences that arise as inner
sums of standard frames of two Hilbert C*-modules with respect to their
properties in the context of the corresponding direct sum Hilbert C*-module.
The main achievement is a description of the various situations in terms of
the relation of the two corresponding orthogonal projections onto the ranges
of the frame transforms in $l_2(A)$. Also, we obtain some uniqueness result
up to similarity for the choice of certain complements in inner sums.
Strong disjointness of standard frames can be described in terms of properties
of their dual frames. Some comments and remarks complete the picture. The main
source of inspiration is chapter two of \cite{HanLarson} where the Hilbert
space situation has been considered, however some changes appear in our more
general setting. We start with the definition of some notions that stand for
the principal situations.

\begin{definition}
  The pairs of standard frames $(\{ x_j \}_j, \{y_j \}_j)$ and $(\{ z_j \}_j,
  \{w_j \}_j)$ of finitely or countably generated Hilbert $A$-modules
  $\mathcal H$ and $\mathcal K$, respectively, are {\it similar} if
  the inner sums $\{x_j \oplus y_j \}_j$ and    \linebreak[4]
  $\{ z_j \oplus w_j \}_j$ are related by an invertible adjointable bounded
  operator       \linebreak[4]
  $T_1 \oplus T_2 :
  \mathcal H \oplus \mathcal K \to \mathcal H \oplus \mathcal K$ by the
  rule $T_1(x_j) \oplus T_2(y_j) = z_j \oplus w_j$ for every $j \in \mathbb J$,
  where $T_1$ acts on $\mathcal H$ and $T_2$ acts on $\mathcal K$.

  Let $\mathcal M$ be a finitely or countably generated Hilbert $A$-module
  that possesses a Riesz basis.
  Let $\{ x_j \}_j$ and $\{ y_j \}_j$ be standard frames of finitely or
  countably generated Hilbert C*-modules $\mathcal H$ and $\mathcal K$,
  respectively, with the property that $\mathcal H \oplus \mathcal K \cong
  \mathcal M$. If both $\{ x_j \}_j$ and $\{ y_j \}_j$ are normalized
  tight and their inner sum $\{ x_j \oplus y_j \}_j$ is a normalized tight
  frame and orthogonal Hilbert basis then the frames $\{ x_j \}_j$ and
  $\{ y_j \}_j$ are {\it strongly complementary to each other}. The two
  frames $\{ x_j \}_j$ and $\{ y_j \}_j$ are {\it strongly complementary to
  each other} if they are similar to a pair of strongly complementary
  normalized tight frames.
  If the inner sum of the frames $\{ x_j \oplus y_j \}_j$ yields a
  Riesz basis of $\mathcal M$ then the frame $\{ y_j \}_j$ is said to be a
  {\it complementary frame to $\{ x_j \}_j$}, and $\{ x_j \}_j$ {\it is
  complementary to $\{ y_j \}_j$}.

  Two standard normalized tight frames $\{ x_j \}_j$ and $\{ y_j \}_j$ of
  $\mathcal H$ and $\mathcal K$, respectively, are {\it strongly disjoint} in
  case the inner sum of them is a standard normalized tight frame of $\mathcal
  M$. Analogously, two standard frames $\{ x_j \}_j$ and $\{ y_j \}_j$ of
  $\mathcal H$ and $\mathcal K$, respectively, are {\it strongly disjoint} if
  they are similar to a strongly disjoint pair of standard normalized
  tight frames of $\mathcal H$ and $\mathcal K$.
  If the inner sum of the frames $\{ x_j \oplus y_j \}_j$ is a standard frame
  of $\mathcal M$ then the frames $\{ x_j \}_j$ and $\{ y_j \}_j$ are said to
  be {\it disjoint}. If the set $\{ x_j \oplus y_j \}_j$ has a trivial
  orthogonal complement in $\mathcal M$ then the frames $\{ x_j \}_j$ and
  $\{ y_j \}_j$ are called {\it weakly disjoint}.
\end{definition}

Examples for the different types of pairs of standard frames can be found
in \cite{HanLarson} for the Hilbert space situation.
All these properties are invariant under unitary equivalence of pairs of frames.
We recall that the property of frames to be standard extends from the summands
to their inner sum. Moreover, by \cite[Prop.~2.2]{LaFr98} any Riesz basis with
frame bounds equal to one is automatically an orthogonal Hilbert basis with
projection-valued inner product values for equal entries taken from the
basis set. If two complementary standard frames $\{ x_j \}_j$ and $\{ y_j \}_j$
of finitely or countably generated Hilbert $A$-modules $\mathcal H$ and
$\mathcal K$, respectively, realize as their inner sum a Riesz basis of
$\mathcal H \oplus \mathcal K$ with frame bounds equal to one then both they
are normalized tight since lower frame bounds of inner summands can only
increase and upper frame bounds only decrease. Therefore, they have to be
strongly complementary. At the contrary, the inner sum of complementary
standard normalized frames need not to yield a Riesz basis with frame bounds
equal to one, in general. To see this consider two complementary standard
normalized tight frames $\{ x_j \}_j$ and $\{ y_j \}_j$ of countably generated
Hilbert $A$-modules $\mathcal H$ and $\mathcal K$, respectively, and suppose
$\mathcal H \oplus \mathcal K \cong l_2(A)$. Then the frame transform $\theta:
\mathcal H \oplus \mathcal K \to l_2(A)$ equals $(\theta_x \oplus 0) + (0 \oplus
\theta_y)$ for the respective frame transforms $\theta_x$ and $\theta_y$ of
the inner summands by construction. Therefore, $\theta \theta^*: l_2(A) \to
l_2(A)$ is the sum of the orthogonal projections $P_{\mathcal H}: l_2(A) \to
\theta_x(\mathcal H)$ and $P_{\mathcal K}: l_2(A) \to \theta_y(\mathcal K)$.
By supposition the sequence $\{ x_j \oplus y_j \}_j$ is a standard Riesz basis
of $l_2(A)$. This forces $\theta \theta^*$ to be invertible, and the ranges
of $P_{\mathcal H}$ and $P_{\mathcal K}$ intersect only in the zero element.
However, for a fixed orthogonal projection $P: l_2(A) \to l_2(A)$ there are
plenty of orthogonal projections $Q: l_2(A) \to l_2(A)$ such that $P(l_2(A))
\cap Q(l_2(A)) = \{ 0 \}$ and $P+Q$ is invertible, not only $Q={\rm id}_{l_2
(A)}-P$, in general. If the frame bounds of the Riesz basis $\{ x_j \oplus y_j
\}_j$ equal to one, then $\theta \theta^*$ is the identity operator on $l_2(A)$
and $P_{\mathcal H} = {\rm id}_{l_2(A)} -P_{\mathcal K}$. Summarizing and
extending our discussion we have the following facts:

\begin{lemma} \label{lemma-6.1} {\rm (cf.~\cite[Lemma 2.4, Prop.~2.5]{HanLarson} ) }
  Let $\{ x_j \}_j$ and $\{ y_j \}_j$ be standard frames of finitely or
  countably generated Hilbert $A$-modules $\mathcal H$ and $\mathcal K$,
  respectively, over a unital C*-algebra $A$.
  If these frames are similar, then the $A$-linear span of the sequence
  $\{ x_j \oplus y_j \}_j$ possesses a non-trivial orthogonal complement in
  the Hilbert $A$-module $\mathcal H \oplus \mathcal K$.
  \newline
  If $\{ x_j \}_j$ and $\{ y_j \}_j$ are normalized tight in addition and the
  sequence $\{ x_j \oplus y_j \}_j$ becomes a standard normalized tight frame
  for the norm-closure of the $A$-linear span of it, then $\overline{{\rm span}
  \{ x_j \oplus y_j \}_j } = \mathcal H \oplus \mathcal K$.
\end{lemma}

\begin{proof}
For the first assertion, consider an invertible adjointable bounded operator
$T: \mathcal H \to \mathcal K$, $T(x_j)=y_j$ for $j \in \mathbb J$.
Note, that the $A$-linear spans of the sequences $\{ y_j \oplus y_j \}_j$ and
$\{ -y_j \oplus y_j \}$ are non-trivial orthogonal to each other subsets of
$\mathcal K \oplus \mathcal K$.
Since $(T \oplus {\rm id})(x_j \oplus y_j) = y_j \oplus y_j$ for all $j$ the
linear span of the set $\{ x_j \oplus y_j \}_j$ cannot be dense in $\mathcal
H \oplus \mathcal K$.

To show the second assertion, fix an index $l \in \mathbb J$. By supposition
 \begin{eqnarray*}
   x_l \oplus y_l & = &
   \sum_{j \in \mathbb J} \langle x_l \oplus y_l, x_j \oplus y_j \rangle
                                          (x_j \oplus y_j) \\
   & = &
   \sum_{j \in \mathbb J} (\langle x_l,x_j \rangle + \langle y_l,y_j \rangle)
                                          (x_j \oplus y_j) \, .
\end{eqnarray*}
Therefore, $0 = \sum_j \langle x_l,x_j \rangle y_j = \sum_j \langle y_l,y_j
\rangle x_j$ by the property of the frames to be normalized tight, and
  \begin{eqnarray*}
    x_l \oplus 0 & = &
    \sum_j \langle x_l \oplus 0, x_j \oplus y_j \rangle (x_j \oplus y_j) \, , \\
    0 \oplus y_l & = &
    \sum_j \langle 0 \oplus y_l, x_j \oplus y_j \rangle (x_j \oplus y_j)
  \end{eqnarray*}
for every $l \in \mathbb J$. This implies the identity
$\overline{{\rm span}\{ x_j \oplus y_j \}_j } = \mathcal H \oplus \mathcal K$.
\end{proof}

We have to point out that for the second statement the supposition to
$\{ x_j \oplus y_j \}_j$ to form a normalized tight frame is essential. The
statement is in general false for arbitrary standard frames, a fact already
emphasized by Deguang Han and D.~R.~Larson in \cite{HanLarson} for the Hilbert
space situation. To give some criteria for the different types of disjointness
in terms of the ranges of the frame transforms we make use of the various
relations of the projections onto them. The considerations use the non-trivial
fact from Hilbert C*-module theory that the equivalence of two norms induced
by two C*-valued inner products on a given Hilbert C*-module induces the
equivalence of the corresponding values of the C*-valued inner products with
equal entries inside the positive cone of the C*-algebra of coefficients,
and thereby the bounds are preserved, cf.~\cite[Prop.~6]{Fr97-1}.

\begin{theorem}  \label{th-disjoint} {\rm (cf.~\cite[Th.~2.9]{HanLarson})}
  Let $\{ x_j \}_j$ and $\{ y_j \}_j$ be standard frames of finitely or
  countably generated Hilbert $A$-modules ${\mathcal H}_1$ and ${\mathcal H}_2$,
  respectively, over a unital C*-algebra $A$. Suppose, $\mathcal H \oplus
  \mathcal K \cong l_2(A)$. Let $P$ and $Q$ the respective projections from
  $l_2(A)$ onto the ranges $\theta_1({\mathcal H}_1)$ and $\theta_2({\mathcal
  H}_2)$ of the corresponding frame transforms. Then:

  \begin{list}{(\roman{marke})}{\usecounter{marke}}
  \item The pair $\{ x_j \}_j$ and $\{ y_j \}_j$ is a strongly complementary
        pair if and only if $P={\rm id}_{l_2(A)}-Q$.
  \item The pair $\{ x_j \}_j$ and $\{ y_j \}_j$ is a complementary pair if and
        only if $P(l_2(A)) \cap Q(l_2(A)) = \{ 0 \}$ and $P(l_2(A)) + Q(l_2(A))
        = l_2(A)$ topologically, i.e.~$(P+Q)$ is invertible.
  \item $\{ x_j \}_j$ and $\{ y_j \}_j$ are strongly disjoint if and only if
        $PQ=QP=0$.
  \item $\{ x_j \}_j$ and $\{ y_j \}_j$ are disjoint if and only if $P(l_2(A))
        \cap Q(l_2(A))$ $= \{ 0 \}$ and $P(l_2(A)) + Q(l_2(A))$ is norm-closed.
  \item $\{ x_j \}_j$ and $\{ y_j \}_j$ are weakly disjoint if and only if
        $P(l_2(A)) \cap Q(l_2(A)) = \{ 0 \}$.
  \end{list}
\end{theorem}

\begin{proof}
By Theorem \ref{th4.2} the frame transforms of similar standard frames of a
fixed Hilbert C*-module have always the same range. So we can restrict our
proof to the consideration of standard normalized tight frames. Note that
the property to be standard transfers from the summands of an inner sum to
the inner sum itself.

We gave already the arguments for the first and the second statements in this
special situation. To demonstrate (iii) and (iv) let $\{ e_j \}_j$ be the
standard orthonormal basis of $l_2(A)$. If $P(l_2(A)) \cap
Q(l_2(A)) = \{ 0 \}$ and the operator $(P+Q)$ has a closed range, then
$0 \leq P+Q \leq 2 \cdot {\rm id}_{l_2(A)}$ by the positivity of $(P+Q)$, and
the canonical norm $\| \langle .,. \rangle_{\mathcal H} + \langle .,.
\rangle_{\mathcal K} \|_A^{1/2}$ on $\theta(\mathcal H \oplus \mathcal K)
\subseteq l_2(A)$ is equivalent to the standard norm on $l_2(A)$. By
\cite[Prop.~6]{Fr97-1} this equivalence of norms forces the equivalence
  \begin{eqnarray} \label{toll}
   C_1 \cdot \langle \theta(x \oplus y), \theta(x \oplus y) \rangle_{l_2(A)}
   & \leq &
   \langle x,x \rangle_{\mathcal H} + \langle y,y \rangle_{\mathcal K}  \\ \nonumber
   & \leq &
   C_2 \cdot \langle \theta(x \oplus y), \theta(x \oplus y) \rangle_{l_2(A)}
  \end{eqnarray}
for some fixed constants $0 < C_1,C_2 < +\infty$ and any $x \in \mathcal H$,
$y \in \mathcal K$. We obtain the following equality:
  \begin{eqnarray*}
    \lefteqn{       \sum_j
      \langle x \oplus y,x_j \oplus y_j \rangle_{\mathcal H \oplus \mathcal K}
      \langle x_j \oplus y_j, x \oplus y \rangle_{\mathcal H \oplus \mathcal K}
    } \\
    & = &
    \sum_j
    |\langle x,x_j \rangle_{\mathcal H} + \langle y,y_j \rangle_{\mathcal K} |^2 \\
    & = &
    \sum_j
       |\langle \theta_x(x),P(e_j) \rangle_{l_2(A)} + \langle \theta_y(y),Q(e_j)
       \rangle_{l_2(A)} |^2 \\
    & = &
    \sum_j
       |\langle \theta_x(x),e_j \rangle_{l_2(A)} + \langle \theta_y(y),e_j
       \rangle_{l_2(A)} |^2 \\
    & = &
    \sum_j |\langle \theta_x(x)+ \theta_y(y),e_j \rangle_{l_2(A)} |^2 \\
    & = &
    \langle \theta_x(x)+\theta_y(y), \theta_x(x)+\theta_y(y) \rangle_{l_2(A)}
    = \langle \theta(x \oplus y) \theta(x \oplus y) \rangle_{l_2(A)}
  \end{eqnarray*}
Together with the equivalence (\ref{toll}) above we get the frame properties
of the inner sum $\{ x_j \oplus y_j \}_j$ for the Hilbert $A$-module
$\mathcal H \oplus \mathcal K$. In case $PQ=QP=0$ the optimal constants
$C_1,C_2$ in (\ref{toll}) are equal to one and the frame $\{ x_j \oplus y_j
\}_j$ is normalized tight. We leave the easy demonstration of the converse
implications in (iii) and (iv) to the reader.

Suppose $\{ x_j \}_j$ and $\{ y_j \}_j$ are weakly disjoint. If there exists
a non-zero element $z \in P(l_2(A)) \cap Q(l_2(A))$ then $z = \theta_x(z_1)
=\theta_y(z_2)$ for some non-zero $z_1 \in \mathcal H$, $z_2 \in \mathcal K$.
The equality
  \begin{eqnarray*}
  \langle z_1 \oplus -z_2, x_j \oplus y_j \rangle_{\mathcal H \oplus \mathcal K}
    & = &
  \langle z_1,x_j \rangle_{\mathcal H} - \rangle z_2,y_j \rangle_{\mathcal K} \\
    & = &
  \langle z,P(e_j) \rangle_{l_2(A)} - \langle z,Q(e_j) \rangle_{l_2(A)} \\
    & = &
  \langle z,e_j \rangle_{l_2(A)} - \langle z,e_j \rangle_{l_2(A)} = 0
  \end{eqnarray*}
$(j \in \mathbb J)$ shows the existence of a non-trivial orthogonal complement
to the $A$-linear span of the sequence $\{ x_j \oplus y_j \}_j$, a contradiction
to the assumptions on it.

Conversely, assume $P(l_2(A)) \cap Q(l_2(A)) = \{ 0 \}$ and the existence
of an element $x \oplus y$ that is orthogonal to the $A$-linear span of
the sequence $\{ x_j \oplus y_j \}_j$. Then
  \begin{eqnarray*}
     0 & = &
      \langle x,x_j \rangle_{\mathcal H} + \langle y,y_j \rangle_{\mathcal K} =
      \langle \theta_x(x),P(e_j) \rangle_{l_2(A)} + \langle \theta_y(y),Q(e_j)
      \rangle_{l_2(A)} \\
      & = &
      \langle \theta(x \oplus y), e_j \rangle_{l_2(A)}   \, ,
  \end{eqnarray*}
and  $\theta(x \oplus y)= \theta_x(x) + \theta_y(y) = 0$. Consequently,
$\theta_x(x) = - \theta_y(y) = 0$, and by the injectivity of frame transforms
we obtain $x=0$ and $y=0$. So the frames $\{ x_j \}_j$ and $\{ y_j \}_j$
are weakly disjoint.
\end{proof}

There is still the possibility that inner sums of weakly disjoint frames
may span either norm-dense subsets or only weakly dense subsets.
Also, in between strong complementarity and complementarity we can consider
the situation that not only $(P+Q)$ is invertible, but also $(2 \cdot
{\rm id}_{l_2(A)} -P-Q) = ({\rm id}_{l_2(A)}-P)+({\rm id}_{l_2(A)}-Q)$ at the
same time. This condition is equivalent to the invertibility of the difference
$(P-Q)$ since $(P-Q)^2=(P+Q)(2 \cdot {\rm id}_{l_2(A)} -P-Q)$ and the latter
two operators always commute, however it does not force $P$ to be equal to
$({\rm id}_{l_2(A)}-Q)$, in general. We will not discuss the details of these
situations here.

The next remarkable fact is some kind of uniqueness of complementary frames.

\begin{proposition} {\rm (cf.~\cite[Prop.~2.1]{HanLarson}) }
  Let $\{ x_j \}_j$ be a standard frame of a finitely or countably generated
  Hilbert $A$-module $\mathcal H$ over a unital C*-algebra $A$. Let $\{ y_j
  \}_j$ and $\{ z_j \}_j$ be standard frames of Hilbert $A$-modules $\mathcal
  K$ and $\mathcal L$, respectively, that are complementary to $\{ x_j \}_j$.
  Suppose, $\langle y_j,y_j \rangle_{{\mathcal M}} = \langle z_j,z_j
  \rangle_{{\mathcal N}}$ for every $j \in {\mathbb J}$.
  Then there exists an invertible operator $T: \mathcal K \to \mathcal L$ such
  that $z_j=T(y_j)$ for every $j \in \mathbb J$, i.e.~these frames are similar.
\end{proposition}

\begin{proof}
Let $T_1,T_2,T_3,T_4$ be invertible adjointable bounded operators such that
the sequences $\{ T_1(x_j) \oplus T_2(y_j) \}_j$ and $\{ T_3(x_j) \oplus
T_4(z_j) \}_j$ are orthogonal bases of the Hilbert $A$-modules $\mathcal H
\oplus \mathcal K$ and $\mathcal H \oplus \mathcal L$, respectively. For the
standard normalized tight frames $\{ T_1(x_j) \}_j$ and $\{ T_3(x_j) \}_j$
we have the identity $T_1(x_j) = T_1T_3^{-1}(T_3(x_j))$ for $j \in \mathbb J$.
By \cite[Prop.~5.10(ii)]{LaFr98} the operator $T_1T_3^{-1}$ has to be unitary.
Consequently, the sequence $\{ T_1(x_j) \oplus T_4(z_j) \}_j$ is an orthogonal
Hilbert basis that is unitarily equivalent to the orthogonal Hilbert basis
$\{ T_3(x_j) \oplus T_4(z_j) \}_j$. By \cite[Prop.~5.2]{LaFr98} there exists a
unitary operator $V: \mathcal K \to \mathcal L$ such that $T_4(z_j)=VT_2(y_j)$
holds for every $j \in \mathbb J$. Finally, we obtain $z_j=T_4^{-1}VT_2(y_j)$
for $j \in \mathbb J$, the desired result.
\end{proof}

\begin{proposition}
  Let $\{ x_j \}_j$ and $\{ y_j \}_j$ be standard frames of the finitely or
  countably generated Hilbert $A$-modules $\mathcal H$ and $\mathcal K$, resp.,
  and let $S_x$ and $S_y$ be their corresponding frame operators. The following
  conditions are equivalent:

  \begin{list}{(\roman{marke})}{\usecounter{marke}}
  \item The frames $\{ x_j \}_j$ and $\{ y_j \}_j$ are strongly disjoint.
  \item The corresponding canonical dual frames $\{ S_x(x_j) \}_j$,
        $\{ S_y(y_j) \}_j$ are strongly disjoint.
  \item The frame operator $S_{x+y}$ of the frame $\{ x_j \oplus y_j \}_j$
        equals the sum of the frame operators $S_x$ and $S_y$, i.e.~the
        canonical dual frame of $\{ x_j \oplus y_j \}_J$ is the inner sum of
        the canonical dual frames of $\{ x_j \}_j$ and $\{ y_j \}_j$.
  \item $\sum_j \langle x,S_x(x_j) \rangle y_j = 0$ for any $x \in \mathcal H$.
  \item $\sum_j \langle y,S_y(y_j) \rangle x_j = 0$ for any $y \in \mathcal K$.
  \item $\sum_j \langle x,x_j \rangle S_y(y_j) = 0$ for any $x \in \mathcal H$.
  \item $\sum_j \langle y,y_j \rangle S_x(x_j) = 0$ for any $y \in \mathcal K$.
  \item $\sum_j \langle x,x_j \rangle y_j = 0$ for any $x \in \mathcal H$.
  \item $\sum_j \langle y,y_j \rangle x_j = 0$ for any $y \in \mathcal K$.
  \end{list}
\end{proposition}

\begin{proof}
By Theorem \ref{th-disjoint} the ranges of the corresponding frame transforms
$\theta_x$ and $\theta_y$ are orthogonal to each other subsets if and only if
condition (i) holds. Since any standard frame is similar to its canonical
dual frame and strong disjointness is invariant under summand-wise similarity
we obtain the equivalence of (i) and (ii).

To show the implication (i)$\to$(iii) observe that the frame transform
$\theta_{x+y}$ of the inner sum $\{ x_j \oplus y_j \}_j$ equals $(\theta_x
\oplus 0_{\mathcal K})+(0_{\mathcal H} \oplus \theta_y)$ since the ranges of
both the frame transforms are orthogonal to each other by Theorem
\ref{th-disjoint},(iii). Consequently, $\theta_{x+y}^*\theta_{x+y} =
\theta_x^*\theta_x \oplus \theta_y^*\theta_y$, and $S_{x+y} = S_x \oplus S_y$.

The frames $\{ S_x(x_j) \}_j$ and
$\{ S_y(y_j) \}_j$ are always standard normalized tight frames of $\mathcal H$
and $\mathcal K$, respectively. If they are strongly disjoint then the inner
sum of them is a standard normalized tight frame of $\mathcal H \oplus
\mathcal K$ since $S = (S_x \oplus 0) + (0 \oplus S_y)$. By Lemma
\ref{lemma-6.1} and its proof the equalities
  \[
     \sum_j \langle x,S_x(x_j) \rangle S_y(y_j) = 0 \quad , \quad
     \sum_j \langle y,S_y(y_j) \rangle S_x(x_j) = 0
  \]
hold in $\mathcal K$ and $\mathcal H$, respectively, for any $x \in \mathcal H$,
$y \in \mathcal K$. Since the operators $S_x$ and $S_y$ are invertible we arrive
at (iv) and (v) omitting them at the right end. We derive (vi) if we
take into account that $\langle S_x(x),x_j \rangle = \langle x,S_x(x_j) \rangle$
for every $x \in \mathcal H$ and $j \in \mathbb J$, where $S_x(x)$ runs over
the entire set $\mathcal H$ if $x$ does so. Analogously we derive (vii), (viii) and
(ix).

Conversely, if one of the formulae (iv)-(ix) holds then the ranges of the
frame transforms of the participating two frames are orthogonal in $l_2(A)$
and they are strongly disjoint by Theorem \ref{th-disjoint}. Again,
summand-wise similarity and invariance of strong disjointness under this
equivalence relation gives (i) and (ii), and therewith the other conditions.

Starting with (iii) and $x \in \mathcal H$, $y \in \mathcal K$ a short
counting shows
  \begin{eqnarray*}
    \lefteqn{ x \oplus y = \sum_j \langle x \oplus y, S_{x+y}(x_j \oplus y_j) \rangle
                                                        (x_j \oplus y_j)} \\
    & = & \sum_j \langle x \oplus y, S_x(x_j) \oplus S_y(y_j) \rangle
                                                        (x_j \oplus y_j) \\
    & = & \sum_j (\langle x,S_x(x_j) \rangle + \langle y,S_y(y_j) \rangle)
                                                        (x_j \oplus y_j) \\
    & = & x \oplus y +
          0 \oplus \left( \sum_j \langle x,S_x(x_j) \rangle  y_j \right) +
          \left( \sum_j \langle y,S_y(y_j) \rangle x_j \right) \oplus 0  \\
    & = & x \oplus y +
          \left( \sum_j \langle y,S_y(y_j) \rangle x_j \right)
          \oplus \left( \sum_j \langle x,S_x(x_j) \rangle  y_j \right)
  \end{eqnarray*}
and, therefore, $0 =  \left( \sum_j \langle y,S_y(y_j) \rangle x_j \right)
\oplus \left( \sum_j \langle x,S_x(x_j) \rangle  y_j \right)$. We get
conditions (iv) and (v) that were already shown to be equivalent to (i).
\end{proof}

We close this section with a proposition that opens a link between special
representations of the Cuntz algebras ${\mathcal O}_n$ on Hilbert C*-modules
$l_2(A)$ and $n$-tuples of pairwise strongly complementary standard normalized
tight frames of $l_2(A)$.

\begin{proposition} {\rm (cf.~\cite[Prop.~2.21]{HanLarson} )}
  Let $\{ x_j \}_j$ and $\{ y_j \}_j$ be standard normalized tight frames
  of the Hilbert $A$-module $l_2(A)$ over a unital C*-algebra $A$. Suppose,
  there are two adjointable bounded operators $T_1,T_2$ on $l_2(A)$ such that
  $T_1T_1^*+T_2T_2^*={\rm id}_{l_2(A)}$. Then the sequence $\{ T_1(x_j)+
  T_2(y_j) \}_j$ is a standard normalized tight frame of $l_2(A)$.
  In particular,
  the sequence $\{ \alpha x_j + \beta y_j \}_j$ is a standard normalized tight
  frame whenever $\alpha, \beta \in \mathbb C$ with $|\alpha|^2+|\beta|^2=1$,
  i.e.~two strongly disjoint standard normalized tight frames are arcwise
  connected within the set of standard normalized tight frames.

  More generally, if $(\{ x_{1j} \}_j, ... ,  \{ x_{nj} \}_j )$ is an $n$-tuple of
  pairwise strong\-ly disjoint standard normalized tight frames of $l_2(A)$ and
  $(T_1,...,T_n)$ is an $n$-tuple of adjointable bounded operators on $l_2(A)$
  with $\sum_j T_jT_j^* = {\rm id}_{l_2(A)}$, then $\{ \sum_{k=1}^n T_k(x_{kj})
  \}_j$ is a normalized tight frame of $l_2(A)$.
\end{proposition}

\begin{proof}
Let $\{ e_j \}_j$ be a standard orthonormal basis of $l_2(A)$, let $\theta_1$,
$\theta_2$ be the frame transforms of the frames $\{ x_j \}_j$ and $\{ y_j \}_j$,
respectively. The operators $\theta_1, \theta_2$ are isometries with orthogonal
ranges in $l_2(A)$ by assumption. Moreover, $\theta_1^*\theta_2=\theta_2^*\theta_1=0$
and
$\theta_1^*(e_j)=x_j$, $\theta_2^*(e_j)=y_j$ for any $j \in \mathbb J$.
Consider the new operator $T=T_1\theta_1^*+T_2\theta_2^*$ and observe that
  \[
     TT^* = T_1\theta_1^*\theta_1T_1^* + T_2\theta_2^*\theta_2T_2^* =
     T_1T_1^*+T_2T_2^* ={\rm id}_{l_2(A)} \, .
  \]
So $T^*$ is an isometry, and $\{ T_1(x_j)+T_2(y_j) = T(e_j) \}_j$ is a
standard normalized tight frame of $l_2(A)$.
\end{proof}

\begin{corollary}  {\rm (cf.~\cite[Prop.~2.16]{HanLarson} )}
   Let $(\{ x_{1j} \}_j,...,\{ x_{nj}\}_j )$ be an $n$-tuple of standard normalized
   tight frames and $\{ e_j \}_j$ be an orthonormal basis of $l_2(A)$.
   If the operators $\{ T_1,...,T_n \}$ are partial isometries on $l_2(A)$
   defined by $T_k(x)= \sum_j \langle x,x_{kj} \rangle e_j$ with pairwise
   orthogonal ranges and $\sum_k T_kT_k^*={\rm id}_{l_2(A)}$,
   then the frame $\{ \sum_{k=1}^n T_k(x_{kj}) \}_j$ is an orthonormal basis
   of $l_2(A)$, i.e.~the frames $(\{ x_{1j} \}_j, ... ,  \{ x_{nj} \}_j )$ are
   strongly complementary. The converse also holds.
\end{corollary}

\begin{proof}
If $T_k(l_2(A)) \bot T_l(l_2(A))$ for any $k \not= l$ then define a new operator
$W$ by
   \[
      W(x_{1j} \oplus ... \oplus x_{nj}) = T_1(x_{1j})+...+T_n(x_{nj})
   \]
for every $j \in \mathbb J$.
Obviously, the equality $\langle W(\sum_k \oplus x_{kj}), W(\sum_k \oplus
x_{kj}) \rangle = \langle \sum_k \oplus x_{kj},\sum_k \oplus x_{kj} \rangle$
holds for every $j \in \mathbb J$ and, hence, $W$ is unitary. Also,
   \[
       W \left( \sum_k \oplus x_{kj})=\sum_k T_k(x_{kj} \right) =
       \sum_k T_kT_k^*(e_j)=e_j
   \]
for every $j \in \mathbb J$. So the st.n.t.~frames
$(\{ x_{1j} \}_j, ... ,  \{ x_{nj} \}_j )$ are strongly complementary.

Conversely, if the st.n.t.~frames $(\{ x_{1j} \}_j, ... ,
\{ x_{nj} \}_j )$ are strongly complementary then
   \[
      \sum_k T_kT_k^*(e_n) = \sum_k T_k(x_{kn}) =
      \sum_j \left(\sum_k \langle x_{kn},x_{kj} \rangle \right) e_j =
      e_n
   \]
for every $n \in \mathbb J$
and, consequently, $\sum_k T_kT_k^* = {\rm id}_{l_2(A)}$, where the operators
$(T_1,...,T_n)$ possess orthogonal ranges by the choice of the frames.
\end{proof}

In fact, for every representation of the Cuntz algebra ${\mathcal O}_n$ by
adjointable bounded operators on the Hilbert C*-module $l_2(A)$ for some unital
C*-algebra $A$ we find an $n$-tuple of strongly complementary standard
normalized tight frames of $l_2(A)$ applying the adjoints of the $n$ partial
isometries generating ${\mathcal O}_n$ to some fixed orthonormal basis of
$l_2(A)$. So these two situations are very closely related.

\section{Frame representations and decompositions}

We establish some module frame representation and decomposition theorems
that are directly derived from operator decomposition properties of elements
of arbitrary unital C*-algebras. The results of this section are inspired by
results by P.~G.~Casazza in \cite{Casazza:1} on the Hilbert space situation,
however resting on our observation that most of the operator decomposition
properties are not only valid in discrete type I von Neumann algebras, but
also in general unital C*-algebras. Obvious changes of the statements on frames
in the more general situation of Hilbert C*-modules are caused by the possible
absence of orthogonal or orthonormal bases in the Hilbert C*-modules under
consideration. During the course of explanation we prove some lemmata that are
of general interest.

\begin{lemma}   \label{Cas-1}
         {\rm   (cf.~\cite[Prop.~2.1]{Casazza:1} )   }
  Let $A$ be a C*-algebra, $\{ {\mathcal H}, \langle .,. \rangle \}$ be a
  Hilbert $A$-module. Every adjointable bounded $A$-linear operator $T \in
  {\rm End}_A^{{*}}({\mathcal H})$ can be decomposed into the real multiple of
  a sum of three unitary operators $U_1,U_2,U_3 \in {\rm End}_A^{{*}}
  ({\mathcal H})$, $T=\lambda (U_1+U_2+U_3)$, with $\lambda > \|T\|$. If $T$
  is surjective and possesses a polar decomposition then it is the arithmetic
  mean of two partial isometries $W_1,W_2$,   \linebreak[4]
  $T= (W_1+W_2)/2$.
\end{lemma}

\begin{proof}
Fix $\varepsilon > 0$. The operator
   \begin{equation}  \label{formula}
     T'=\frac{{\rm id}_{{\mathcal H}}}{2} + \frac{1}{2(1+\varepsilon)\|T\|} T
   \end{equation}
is adjointable and invertible since
   \[
     \| {\rm id}_{{\mathcal H}} - T'\| =
     \left\| \frac{{\rm id}_{{\mathcal H}}}{2} -
              \frac{1}{2(1+\varepsilon)\|T\|} T \right\| \leq
     \frac{1}{2} + \frac{1}{2(1+\varepsilon)}  < 1  \, .
   \]
Because of the invertibility of $T'$ in the C*-algebra ${\rm End}_A^{*}
({\mathcal H})$ the operator possesses an polar decomposition $T'=V |T'|$
where $V$ is a unitary operator. Since $\| T' \| < 1$ we obtain
$\| \, |T'| \, \| < 1$ and
   \begin{equation}  \label{formula2}
     | T' | = \frac{1}{2} \left(
        ( |T'| + {\bf i} \cdot | {\rm id}_{{\mathcal H}} - |T'|^2  |^{1/2}) +
        ( |T'| - {\bf i} \cdot | {\rm id}_{{\mathcal H}} - |T'|^2  |^{1/2})
        \right)
   \end{equation}
where $U =  (|T'| + {\bf i} \cdot | {\rm id}_{{\mathcal H}} - |T'|^2  |^{1/2})$
is unitary.
Consequently, $T = (1+\varepsilon) \| T \| (VU+VU^*- {\rm id}_{{\mathcal H}})$
is the sought decomposition by (\ref{formula}).

\noindent
If $T$ is a adjointable bounded $A$-linear map of $\mathcal H$ onto itself that
has a polar decomposition $T = V |T|$ inside ${\rm End}_A^{*}({\mathcal H})$,
then $V$ is a partial isometry with range $\mathcal H$ and ${\rm ker}(V)
= {\rm ker}(|T|)$. Replacing $T'$ by $T$ in (\ref{formula2}) we get $T = 1/2
\cdot (VU+VU^*)$ with partial isometries $VU,VU^*$.
\end{proof}

P.~G.~Casazza gave an example of a normalized tight frame of a separable
Hilbert space that cannot be written as any linear combination of two
orthonormal sequences of that Hilbert space, \cite[Ex.~2.3]{Casazza:1}.
Since $l_2(A) \cong \overline{A \odot l_2}$ for every C*-algebra $A$ there
exists a standard normalized tight frame of $\mathcal H = l_2(A)$ with the
same property. The more remarkable is the following decomposition property
of standard frames of finitely or countably generated Hilbert C*-modules:

\begin{theorem} {\rm   (cf.~\cite[Prop.~2.7]{Casazza:1})  }
  Let $A$ be a unital C*-algebra and $\mathcal H$ be a finitely or countably
  generated Hilbert $A$-module. Every standard frame of $\mathcal H$ is a
  multiple of the sum of two standard normalized tight frames of $\mathcal H$.
\end{theorem}

\begin{proof}
Consider the adjoint $\theta^*: l_2(A) \to \mathcal H$ of the frame operator
$\theta$ induced by the frame $\{ f_j \}$ of $\mathcal H$. It possesses a
polar decomposition and is surjective. Consequently, $\theta^*$ is a linear
combination of two partial isometries by the previous lemma. Since $\theta^*$
maps the elements of the standard orthonormal basis precisely to the frame
elements, and since partial isometries applied to orthonormal bases yield
standard normalized tight frames of the image module by Theorem \ref{th3}
we get the desired result.
\end{proof}

\begin{lemma}  \label{Cas-2} {\rm   (cf.~\cite{Casazza:1,KPR})  }
  Let $A$ be a C*-algebra and $\{ {\mathcal H}, \langle .,. \rangle \}$ be a
  Hilbert $A$-module. The subset of all open (resp., surjective) [adjointable]
  bounded $A$-linear operators on $\mathcal H$ forms an open subset of
  ${\rm End}_A^{[*]}({\mathcal H})$ with respect to the operator norm topology.
  The subset of all invertible [adjointable] bounded $A$-linear operators forms
  a clopen (i.e.~both open and closed) subset of the set of all surjective
  [adjointable] bounded $A$-linear operators on $\mathcal H$.
\end{lemma}

\begin{proof}
We can replace the open maps by the surjective ones in
\cite[Prop.~7.8]{KPR} without changes in the proof, a fact first observed by
P.~G.~Casazza in \cite{Casazza:1}. Hence, we know that the open (resp.,
surjective) bounded linear mappings on $\mathcal H$ form an open subset of the
set of all bounded linear mappings on $\mathcal H$. Intersecting the complement
of this open set with the closed subset of all [adjointable] bounded $A$-linear
operators on $\mathcal H$ and taking again the complement we obtain that the
open (resp., surjective) [adjointable] bounded $A$-linear operators on
$\mathcal H$ form an open subset in ${\rm End}_A^{[*]}({\mathcal H})$.

Since the subset of all invertible [adjointable] bounded $A$-linear operators
is an open subset of ${\rm End}_A^{[*]}({\mathcal H})$ by the spectral theory of
Banach algebras we conclude that the invertible [adjointable] bounded
$A$-linear operators form an open subset of the subset of all surjective
[adjointable] bounded $A$-linear mappings on $\mathcal H$. Again, by
\cite[Prop.~7.9]{KPR} the set of all invertible bounded linear mappings on
$\mathcal H$ becomes a closed subset of the subset of all surjective bounded
linear operators on $\mathcal H$. Intersecting with the closed in
${\rm End}_{{\bf C}}({\mathcal H})$ subset ${\rm End}_A^{[*]}({\mathcal H})$
the analogous statement for [adjointable] bounded $A$-linear maps yields.
\end{proof}

\begin{lemma}  \label{Cas-3}
         {\rm (cf.~\cite[Prop.~2.4]{Casazza:1} )}
  Let $A$ be a C*-algebra, $\{ {\mathcal H}, \langle .,. \rangle \}$ be a
  Hilbert $A$-module. A surjective bounded $A$-linear operator $T : {\mathcal
  H} \to {\mathcal H}$ can be written as the linear combination of two unitary
  operators if and only if $T$ is invertible and adjointable.
\end{lemma}

\begin{proof}
If $T$ is adjointable and invertible then $T=1/2 \cdot (U_1+U_2)$ by the
second part of Lemma \ref{Cas-1}, where $U_1,U_2 \in {\rm End}_A^{*}
({{\mathcal H}})$ are unitary operators.
Conversely, suppose  $T=\lambda U_1 + \mu U_2$ for complex numbers $\lambda,\mu$
and unitary operators $U_1,U_2 \in {\rm End}_A^*({{\mathcal H}})$. Then $T$ is
naturally adjointable. If either $\mu$ or $\lambda$ are zero, then $T$ is
invertible.
So, without loss of generality we start the investigation with $T=U_1+\nu U_2$
for $|\nu| \geq 1$ and unitary operators $U_1,U_2$ dividing our original
decomposition by that complex number of $\{ \lambda, \mu \}$ with the smaller
absolute value. Set $T_t:=tU_1 + (1-t)\nu U_2$ for real $0 \leq t < 1/2$.
Since
    \[
      \|T_t(x)\| \geq (1-t) |\nu| \|U_2(x)\|-t\|U_1(x)\| =((1-t)|\nu|-t)\| x \|
    \]
for every $x \in \mathcal H$ and since $((1-t)|\nu|-t) > 0$ for $0 \leq t < 1/2$
we obtain the injectivity of $T_t$ for every $t \in [0,1/2)$. On the other
hand, $2 \cdot T_{1/2}=T$ and $T$ was supposed to be surjective. Applying
Lemma \ref{Cas-2} we obtain for $t$ close enough to $1/2$ that the operator
$T_t$ has to be surjective since the sequence $\{ T_t \}$ converges to $T$ in
norm for $t \to 1/2$. Furthermore, $T=2 \cdot T_{1/2}$ has to be invertible.
\end{proof}

Now we are in the position to characterize all standard frames of the Hilbert
C*-module $l_2(A)$ that can be written as linear combinations of
two orthonormal bases. For a more general result see Theorem \ref{Cas-98}
below. Note that the proof given works equally well for frames of $A^n$,
$n \in \mathbb N$, that contain exactly $n$ elements. However, this is an
even more exceptional situation, and we prefer to give only the formulation
for the infinitely generated case.

\begin{theorem} \label{Cas-99} {\rm (cf.~\cite[Prop.~2.5]{Casazza:1}) }
   A standard frame $\{ f_j \}$ of a Hilbert $A$-module $\mathcal H = l_2(A)$
   can be written as a linear combination of two orthonormal bases of
   $\mathcal H$ if and only if it is a standard Riesz basis.
\end{theorem}

\begin{proof}
Consider the frame transform $\theta: \mathcal H \to \mathcal H$ of this frame.
Then by the previous result its adjoint $\theta^*$ is invertible if and only
if it is a linear combination of two unitary operators. However, every
orthonormal basis of $\mathcal H$ is connected to the fixed standard orthonormal
basis of $\mathcal H$ by a unitary operator. What is more, $\theta^*(e_j) = f_j$
for every $j \in J$ and the standard orthonormal basis $\{ e_j \}$ of $\mathcal H$.
So the claimed decomposition property of the frame $\{ f_j \}$ is equivalent
to the invertibility of $\theta^*$, and $\theta$ has to be invertible, too.
This in turn shows the property of the frame $\{ f_j \}$ to be a basis.
\end{proof}

Since for countably generated Hilbert C*-modules the adjointability of all
bounded module operators between them is far from guaranteed, in general
(cf.~\cite[Th.~4.3]{Fr1}), and Riesz bases may not exist by Example \ref{example1} we
can show the following result only under restrictive circumstances.

\begin{proposition} {\rm    (cf.~\cite[Cor.~2.2]{Casazza:1}) }
   Let $A$ be a unital C*-algebra and $\mathcal H = l_2(A)$.
   If $\{ x_j : j \in {\bf J} \}$ is a standard frame for $\mathcal H$ with
   upper frame bound $D$ then for every $\varepsilon > 0$ there exist three
   orthonormal bases $\{ f_j \}$, $\{ g_j \}$, $\{ h_j \}$ of $\mathcal H$
   such that $x_j =  D(1+\varepsilon)(f_j+g_j+h_j)$ for every $j \in {\bf J}$.
\end{proposition}

\begin{proof}
Let $\{ e_j \}$ be an orthonormal basis of $\mathcal H$. The operator $\theta^*$
defined by the formula $\theta^*(e_j)=x_j$ is the adjoint of the frame
transform $\theta$, $\theta^*$ is surjective and $\| \theta^* \|$ equals the
upper frame bound $D$ of the frame $\{ x_j \}$. It can be represented as
$\theta^* = D (1+ \varepsilon) (U_1+U_2+U_3)$ for three unitary operators
$U_1,U_2,U_3$ on $\mathcal H$ that depend on the choice of $\varepsilon > 0$,
see Lemma \ref{Cas-1}. Setting $f_j=U_1(e_j)$, $g_j=U_2(e_j)$, $h_j=U_3(e_j)$
we are done.
\end{proof}

Better results can be obtained if we dilate the given Hilbert C*-module into
a bigger one and consider linear combinations of bases and frames therein.

\begin{proposition}  {\rm (cf.~\cite[Prop.~2.8]{Casazza:1}) }
   Let $A$ be a unital C*-algebra and $\mathcal H$ be a finitely or countably
   generated Hilbert $A$-module. For every standard normalized tight frame
   $\{ h_j \}$ of $\mathcal H$ the symbol $\theta$ denotes the corresponding
   frame transform that realizes an isometric embedding of $\mathcal H$ into
   $l_2(A)$.
   Then there exist two orthonormal bases $\{ f_j \}$, $\{ g_j \}$ of $l_2(A)$
   such that $\theta(h_j) = 1/2 \cdot (f_j+g_j)$ for $j \in \mathbb N$.
\end{proposition}

\begin{proof}
Denote by $P$ the (existing) projection of $l_2(A)$ onto the range of $\theta$.
Note, that $\theta(h_j)=e_j$ for $j \in \mathbb N$ and for the fixed orthonormal
basis $\{ e_j \}$ of $l_2(A)$ used for the definition of the frame transform.
Then $f_j:= P(e_j)+({\rm id}_{l_2(A)}-P)(e_j)=e_j$ and $g_j:=P(e_j)-
({\rm id}_{l_2(A)}-P)(e_j)$ are both orthonormal bases of $l_2(A)$, and
$\theta(h_j)=f_j+g_j$ for every $j \in \mathbb N$.
\end{proof}

Every standard frame is the image of a standard normalized tight frame
under an adjointable invertible operator $T$ on the fixed Hilbert $A$-module
$\mathcal H$ by \cite{LaFr98}. Dilating $T$ to $\left( \theta T \theta^{-1} P+
( {\rm id}_{l_2(A)}-P) \right)$ on $l_2(A)$ we get another corollary.

\begin{corollary}
   Let $A$ be a unital C*-algebra and $\mathcal H$ be a finitely or countably
   generated Hilbert $A$-module. For every standard frame
   $\{ h_j \}$ of $\mathcal H$ the symbol $\theta$ denotes the corresponding
   frame transform that realizes an embedding of $\mathcal H$ into $l_2(A)$.
   Then there exist two Riesz bases $\{ f_j \}$, $\{ g_j \}$ of $l_2(A)$
   such that $\theta(h_j) = 1/2 \cdot (f_j+g_j)$ for $j \in \mathbb N$.
\end{corollary}

We finish our considerations with an extension of Theorem \ref{Cas-99} that
has been motivated by the situation for Hilbert spaces described in
\cite{Casazza:1}.

\begin{theorem}  \label{Cas-98} {\rm   (cf.~\cite[Prop.~2.10]{Casazza:1} ) }
   Let $A$ be a unital C*-algebra. Every standard frame of $\mathcal H =
   l_2(A)$ is a multiple of the sum of an orthonormal basis of $\mathcal H$
   and a Riesz basis of $\mathcal H$.
\end{theorem}

\begin{proof}
The proof is given in the same manner as for Lemma \ref{Cas-1}, with slight
modifications.
If $\theta^*: l_2(A) \to l_2(A)$ is the adjoint operator of the frame transform
arising from the given standard frame $\{ f_j \}$ of $l_2(A)$ then define
an operator $T'$ by the formula $4T' = {\rm id}_{l_2(A)} + (1-\varepsilon) T /
\| T \|$ for any fixed $\varepsilon > 0$. Again, $\| {\rm id}_{l_2(A)} -T' \|
< 1$ and $\|T'\| \leq 1$ forces $T'$ to be invertible and to possess a
representation $T' = 1/2 \cdot (W+W^*)$ for a unitary operator $W$ on $l_2(A)$
by (\ref{formula2}). (Note that $W$ replaces $VU$ of the proof of Lemma
\ref{Cas-1}.) Resolving the equation for $T'$ with respect to $T$ we obtain
   \[
      \theta^* = \frac{2 \| \theta^* \|}{(1-\varepsilon)} \left( W + W^* -
      \frac{3}{2} {\rm id}_{l_2(A)} \right) \, .
   \]
Since $W$ is unitary the sequence $\{ W(e_j) \}$ forms an orthonormal basis
of $l_2(A)$, where $\{ e_j \}$ denotes the standard orthonormal basis of
$l_2(A)$ used to build the frame transform $\theta$. The second summand
$(W^*-3/2 \cdot {\rm id}_{l_2(A)})$ is injective and adjointable. Since
   \[
      \left\|
       {\rm id}_{l_2(A)} - \left( \frac{-1}{2} \left( W^*-\frac{3}{2}
       {\rm id}_{l_2(A)} \right) \right)
      \right\| =
      \left\| \frac{1}{4} {\rm id}_{l_2(A)} + \frac{1}{2} W^* \right\| < 1
   \]
the operator $(W^*-3/2 \cdot {\rm id}_{l_2(A)})$ multiplied by $-1/2$ has to be
invertible, i.e.~it is also surjective.
Consequently, the sequence $\{ (W^*-3/2 \cdot {\rm id}_{l_2(A)})(e_j) \}$
forms a standard Riesz basis of $l_2(A)$. In total we obtain
   \[
      f_j = \theta^*(e_j) = \frac{2 \| \theta^* \|}{(1-\varepsilon)} \left(
      W(e_j) + (W^* - 3/2 \cdot {\rm id}_{l_2(A)})(e_j) \right) \quad ,
      \quad j \in \mathbb N \, ,
   \]
as desired.
\end{proof}

Again, the proof works equally well for frames of Hilbert $A$-modules $A^n$,
$n \in \mathbb N$, consisting exactly of $n$ elements, a rather exceptional
situation.



\end{document}